\documentclass[letterpaper, 10 pt, conference]{ieeeconf}  

\IEEEoverridecommandlockouts                              

\usepackage{amsmath} 
\usepackage{amssymb}  
\usepackage{xcolor} 
\usepackage{float}
\usepackage{graphicx}

\usepackage{tikz}
\usepackage{textcomp}
\usepackage{hyperref}
\usepackage{lipsum}

\usepackage{amsthm}
\theoremstyle{definition}
\newtheorem{definition}{Definition}[section]

\newcommand{\R}{\mathbb{R}}

\newcommand\copyrighttext{%
  \footnotesize \textcopyright 2023 IEEE. Personal use of this material is
  permitted. Permission from IEEE must be obtained for all other uses, in any
  current or future media, including reprinting/republishing this material for
  advertising or promotional purposes, creating new collective works, for resale
  or redistribution to servers or lists, or reuse of any copyrighted component
  of this work in other works.
}

\newcommand\copyrightnotice{%
  \begin{tikzpicture}[remember picture,overlay]
  \node[anchor=south,yshift=10pt] at (current page.south) {\fbox{\parbox{\dimexpr\textwidth-\fboxsep-\fboxrule\relax}{\copyrighttext}}};
  \end{tikzpicture}%
}

\title{\LARGE \bf
  Topological Data Analysis of Electroencephalogram Signals for Pediatric
  Obstructive Sleep Apnea
}

\author{Shashank Manjunath$^{1}$, Jose A. Perea$^{2}$ and Aarti Sathyanarayana$^{3}$
  \thanks{*J. A. Perea was partially supported by the National Science
  Foundation through grants CCF-2006661 and CAREER award DMS-1943758.}
	\thanks{$^{1}$Shashank Manjunath is with the Khoury College of Computer
		Sciences, Northeastern University, Boston, MA, 02115
			{\tt\small manjunath.sh@northeastern.edu}}%
  \thanks{$^{2}$Jose A. Perea is with the Department of Mathematics, Northeastern University, Boston, MA, 02115,
    and the Khoury College of Computer Sciences, Northeastern University, Boston, MA, 02115
			{\tt\small j.pereabenitez@northeastern.edu}}%
	\thanks{$^{3}$Aarti Sathyanarayana is with the Khoury College of Computer
		Sciences, Northeastern University, Boston, MA, 02115, the Bouvé College of
    Health Sciences, Northeastern University, Boston, MA, 02115, and the Department of
		Biostatistics, Harvard School of Public Health, Boston, MA, 02115
			{\tt\small a.sathyanarayana@northeastern.edu}}%
}

\begin{document}

\maketitle
\copyrightnotice

\thispagestyle{empty}
\pagestyle{empty}

\begin{abstract}

  Topological data analysis (TDA) is an emerging technique for biological signal
  processing. TDA leverages the invariant topological features of signals in a
  metric space for robust analysis of signals even in the presence of noise. In
  this paper, we leverage TDA on brain connectivity networks derived from
  electroencephalogram (EEG) signals to identify statistical differences between
  pediatric patients with obstructive sleep apnea (OSA) and pediatric patients
  without OSA. We leverage a large corpus of data, and show that TDA enables us
  to see a statistical difference between the brain dynamics of the two groups.
  \newline

  \indent \textit{Clinical relevance}— This establishes the potential of
  topological data analysis as a tool to identify obstructive sleep apnea
  without requiring a full polysomnogram study, and provides an initial
  investigation towards easier and more scalable obstructive sleep apnea
  diagnosis.
\end{abstract}

\section{INTRODUCTION}\label{sec:intro}

Recent advancements in signal processing technology, including topological data
analysis (TDA), provide a powerful method for analysis of EEG signals.
TDA leverages ideas from the mathematical field of
topology, and applies these ideas to the analysis of real-world
signals\cite{chazalIntroductionTopologicalData2021b}. Broadly, TDA allows us to
exploit the topological and geometric structures inherent to data, and use these
structures to study fundamental differences between the EEG signals of
obstructive sleep apnea (OSA) positive and OSA negative patients. In this work,
we present TDA techniques which allows identification of OSA using only EEG
signals.

The key assumption in this work is that the brain connectivity network for OSA
positive patients has a fundamentally different topological structure than the
brain connectivity network for OSA negative patients.

Existing techniques for the identification of OSA in children involve a full
overnight sleep study called a polysomnogram (PSG). This requires the patient to
either go to a sleep lab in a medical facility or schedule an overnight sleep
test at home, both of which can take months to schedule.

A successful application of TDA to differentiating OSA positive and OSA negative
patients directly through EEG signals, as we present in this paper, is a strong
step towards easier and less invasive diagnosis of OSA in children, leading to
better patient outcomes.

\section{BACKGROUND}\label{sec:bg}

\subsection{Obstructive Sleep Apnea}\label{subsec:osa}

Pediatric Obstructive Sleep Apnea (OSA) is difficult to diagnose and leads to
significant complications, including poor rest and heart
disease\cite{capdevilaPediatricObstructiveSleep2008a}. In children especially,
OSA can lead to behavioral changes and further reduction in quality of life.
Identification of OSA in patients typically requires an overnight PSG
study~\cite{schechterTechnicalReportDiagnosis2002a}. PSG studies include a
number of sensors, such as EEG sensors, pulse oximetry, and respiratory rate
sensors, among others. Acquisition of this information is difficult for the
patient, as they must sleep in an unknown environment while being observed by
medical personnel, as well as sleep with multiple sensors connected to their
body. PSG studies are considered the gold standard for OSA diagnosis, and
attempts to use other clinical methods such as questionnaires lead to poor
diagnostic accuracy\cite{schechterTechnicalReportDiagnosis2002a}.

Snoring is relatively common in children, with about 10\% incidence in children
ages 1-9. Obstructive sleep apnea occurs in 2-3\% of children with
habitual snoring~\cite{capdevilaPediatricObstructiveSleep2008a}. Further
complication is caused by recent increases in pediatric obesity. Pediatric OSA
is typically caused by adenotonsillary hypertrophy, and has been
correlated with poor growth and failure to
thrive~\cite{chanObstructiveSleepApnea2004a}. Furthermore, daytime drowsiness
(somnolescence) is much less common in children with OSA as compared to adults
with OSA, leading to further difficulties with the diagnosis of OSA.

\subsection{Topological Data Analysis}\label{subsec:tda}

TDA aims to identify quantitative information about the structure of data, and
leverage this information for downstream data analysis
tasks~\cite{chazalIntroductionTopologicalData2021b}. Data is assumed to exist in
a metric space, allowing us to compute distances between different observations,
which are typically viewed as points in an ambient space --- i.e., forming a
point cloud. From the data, one builds a shape on the input point cloud by
varying the distance with which we consider two points as connected (further
discussed in Section~\ref{sec:tda_methods}).

TDA has previously been used to understand the underlying dynamics of the brain
in work
by~\cite{bouraknaTopologicalDataAnalysis2022a,lellaPersistentHomologyFractal2021}.
However, these papers focus on identification of attention-deficit/hyperactivity
disorder (ADHD) in the case of Bourakna et.\
al.~\cite{bouraknaTopologicalDataAnalysis2022a}, and the classification of
specific EEG signal morphology during sleep in the case of Lella et.\
al.~\cite{lellaPersistentHomologyFractal2021}. Work by Lella et.\ al.\ shows the
strength of TDA as an analysis technique for single EEG signals, as this work
outperforms previous state-of-the-art classification techniques with greater
robustness to inherent noise than previous work; however, this work is focused
on the specific morphology of EEG signals. While specific EEG morphology is
useful, in the case of OSA we require a deeper understanding of brain dynamics
in order to appropriately identify OSA-specific features. Additionally, we
require a coherent analysis, involving all collected EEG channels, in order to
fully exploit EEG data for TDA analysis. Work by Bourakna et.\ al.\ provides a
strong foundation for leveraging TDA for analysis of ADHD data, while also using
TDA to understand underlying brain network dynamics. This work performed a
coherent analysis of multiple EEG channels in order to understand brain dynamics
of ADHD patients rather than specific morphologies of individual EEG signals,
and is therefore more robust to noise. We aim to apply these ideas to the
identification of OSA.

OSA is not currently diagnosed by EEG signals in isolation, but rather in
combination with other sensors collected during PSG studies. In this work, we
leverage TDA applied to EEG signals without incorporating other sensors in order
to identify statistically significant differences between OSA positive and OSA
negative groups. This work provides a stepping stone towards improving
scalability and ease of OSA diagnosis for both the clinician and the patient.

\section{METHODS}\label{sec:methods}

\subsection{Data}\label{sec:data}

In this work, we leverage the Nationwide Children's Hospital (NCH) Sleep
DataBank
(SDB)\cite{zhangNationalSleepResearch2018a,leeLargeCollectionRealworld2022a}.
This dataset includes 3,984 PSG studies from 3,673 pediatric patients. This
dataset contains a significant volume of data from a number of sensors,
including EEG, electromyogram, electrooculogram, electrocardiogram, nasal and
oral sensors to measure airflow, and pulse oximetry, among others. This large
volume of data and sensors allows us to study pediatric OSA across a wide
variety of pediatric patients. The dataset additionally contains annotations of
sleep events, such as NREM1 sleep, NREM2 sleep, NREM3 sleep, and REM sleep, in
30 second intervals. We focus on EEG signals collected at 256 Hz using 7-channel
EEG studies, which restricts us to 2,883 sleep studies. We identify patients as
being in the ``Apnea'' or ``No Apnea'' based on the presence of an Apnea event
string annotation in the patient's respective polysomnogram in the dataset. If
an apnea-related string annotation is in the patient's polysomnogram, at any
point during the sleep study, the whole study is placed into the ``Apnea''
group. If no apnea-related string annotation is found, the study is placed into
the ``No Apnea'' group.

\subsection{Signal Preprocessing}\label{sec:signal_proc_methods}

First, we perform basic signal processing on our EEG signals. We filter out 60Hz
and 120Hz power noise using a 3rd order Butterworth bandstop filter.

We then leverage techniques from Bourakna et al. to create distance matrices
from our EEG data\cite{bouraknaTopologicalDataAnalysis2022a}. We process data on
a per-study basis. For a single PSG study we first extract all available EEG
channels, leading to a set of time-series signals $X_i (t)$ for
channel $i \in V$ and time $t \in \{1, \cdots, T\}$. We break these signals into
30-second chunks to correlate with labels of Awake, NREM1, NREM2, NREM3, and REM sleep.
We then create a smoothed periodogram from each channel for each
30-second chunk. A smoothed periodogram allows us to understand the frequency
spectrum of our signal, while smoothing the frequency spectrum slightly using a
small smoothing kernel to eliminate noise.

To calculate the smoothed periodogram, we use the following equations:

\begin{equation}\label{eq:fourier_coef}
	d(\omega_k) = \frac{1}{\sqrt{T}}\sum_{t=1}^T X_i (t) \exp(-i 2 \pi t \omega_k)
\end{equation}

\begin{equation}\label{eq:smooth_periodogram}
	\hat f(\omega_k) = \sum_{\omega} k_h (\omega - \omega_k) d(\omega_k)
	d(\omega_k)^*
\end{equation}

\noindent $k_h (\omega - \omega_k)$ indicates a smoothing kernel centered around
$\omega_k$ and $h$ is a bandwidth parameter. From these smoothed periodograms,
we calculate the squared signal coherence between EEG channels $i$ and $j$:

\begin{equation}\label{eq:cross_correlation}
	\mathcal{C}(X_i, X_j, \omega) =
	\frac{|f_{i,j}(\omega)|^2}{f_{i,i}(\omega)f_{j,j}(\omega)}
\end{equation}

\noindent We then leverage a decreasing function to create a distance between coherence
values for a given frequency $\omega$:

\begin{equation}\label{eq:distance_matrix}
	\mathcal{D}(X_i, X_j, \omega) = 1 - \mathcal{C}(X_i, X_j, \omega)
\end{equation}

The matrix $\mathcal{D}$ is indexed by both channel $i$, channel $j$, and
frequency $\omega$. EEG data is clinically analyzed using the following
frequency bands: Delta band from 0.5 Hz to 4 Hz, Theta band from 4 Hz to 8 Hz,
Alpha band from 8 Hz to 12 Hz, Beta band from 12 Hz to 30 Hz, and Gamma band
from 30 Hz to 50 Hz.

We average across frequency all matrices which occur within a frequency band in
order to obtain a single $7 \times 7$ matrix for that frequency band. Therefore,
for a single 30-second EEG signal from 7 channels, we obtain 5 separate distance
matrices, one for each frequency band. We treat distance matrices separately for
each frequency band and sleep state. These distance matrices serve as the
input for our TDA methods.

\subsection{Topological Data Analysis Methods}\label{sec:tda_methods}

In order to identify if the brain connectivity networks are fundamentally
topologically different for OSA positive and OSA negative patients, we must
understand how to differentiate topological spaces. We first begin with the
definition of a simplicial complex.

\begin{definition}[Simplicial Complex~\cite{edelsbrunnerComputationalTopologyIntroduction2010}]
    A simplicial complex is a collection $K$ of nonempty finite sets $\sigma \in
    K$, called simplices, such that $\sigma \in K$ and $\emptyset \neq \tau
    \subseteq \sigma$ always imply that $\tau \in K$. If $\sigma \in K$ has
    $n+1$ elements (i.e., $\sigma = \{ x_0, \ldots, x_n\}$) then we call it an
    $n$-simplex.
\end{definition}

Intuitively, a 0-simplex $\{x_0\}$ represents a vertex, a 1-simplex
$\{x_0,x_1\}$   an edge between $x_0$ and $x_1$, a 2-simplex $\{x_0,x_1,x_2\}$
encodes a triangle spanned by these vertices, and so on in higher dimensions.
Homology, which we define next, is a descriptor of   the topology of  simplicial
complexes. It measures properties such as the number of connected components,
the existence of holes, voids, and their higher dimensional analogs. Indeed,
fix a field $\mathbb{F}$ (e.g., the rationals or the integers modulo a prime)
and let $C_n(K)$ denote the $\mathbb{F}$-vector space generated by the
$n$-simplices of $K$. That is, each $\gamma \in C_n(K)$ can be written  uniquely
as $\gamma = \sum\limits_{j=0}^{J} a_j \sigma_j$ for scalars $a_j \in
\mathbb{F}$ and $n$-simplices  $\sigma_j \in K$. Define the boundary of an
$n$-simplex $\sigma = \{x_0, \ldots ,x_n\}\in K$ as $\partial_n (\sigma) =
\sum\limits_{i=0}^n (-1)^i \left(\sigma \smallsetminus \{x_i\}\right) \in
C_{n-1}(K)$. Since the $n$-simplices of $K$ form a basis for $C_n(K)$, then we
obtain   linear transformations $\partial_n : C_n(K) \longrightarrow C_{n-1}(K)$
which, as one can check, satisfy $\partial_{n}\circ \partial_{n+1}  = 0 $ for
every $n\in \mathbb{N}$.

\begin{definition}[Simplicial Homology Group~\cite{chazalIntroductionTopologicalData2021b}]
  The $n$th homology group of a simplicial complex $K$ is
  the quotient vector space

  \[
    H_n(K) = \mathsf{Ker}(\partial_n) / \mathsf{Img}(\partial_{n+1})
  \]
\end{definition}

The homology of $K$ can be interpreted as follows:
the dimension of $H_0(K)$ is exactly the number of connected components of $K$,
$\mathsf{dim}(H_1(K))$ counts the essentially distinct closed loops in $K$ bounding an empty region,
$\mathsf{dim}(H_2(K))$ counts cavities,
and similarly for $\mathsf{dim}(H_n(K))$, $n\geq 1$.

In this work, we attempt to understand the topological features of each homology
group (e.g., connected components, holes, cavities, etc.) for the brain
connectivity network of OSA positive and OSA negative patients. Further
treatment of topological concepts can be found
in~\cite{edelsbrunnerComputationalTopologyIntroduction2010}.

In order to calculate homological features on sets of discrete data points, we
apply persistent homology to the distance matrices described in
Section~\ref{sec:signal_proc_methods}.
Indeed, once we have obtained a distance matrix for a
30-second signal, we build a  simplicial complex, in
particular a \emph{ Rips Complex},   and furthermore the
\emph{Rips filtration} using the Python package
\texttt{Ripser.py}\cite{tralieRipserPyLean2018}.

Let $\mathcal{S}$ be our point cloud, with  pairwise distances described by our
previous distance matrix $\mathcal{D}$. The Rips complex $\mathcal{R}$ of
$\mathcal{S}$ and distance $r$ is given by the  simplicial complex
consisting of all subsets of $\mathcal{S}$ with diameter at most $r$, i.e.

\[
  \mathcal{R}(\mathcal{S}, r) = \{\sigma \subset \mathcal{S} | \text{diam}(\sigma) \leq r\}
\]

Note that $\mathcal{R}$ is inherently tied to a given $r$ value. We construct
a \emph{filtration}, called a Rips filtration, which is
the set of Rips complexes created by varying the free parameter $r$. We start
with a small $r$ ($r = 0$), and increase $r$ continuously to obtain a set
of subcomplexes of our simplicial complex\cite{edelsbrunnerComputationalTopologyIntroduction2010}:

\[
  \emptyset = K_0 \subseteq K_1 \subseteq K_2 \subseteq \cdots \subseteq K_m
\]

We leverage \texttt{Ripser.py} to compute our Rips filtration. As we vary $r$
from $K_{i-1}$ to $K_i$, our topological features change, and we gain or lose
topological features\cite{pereaTopologicalTimeSeries2019a}. We quantify the $r$
at which a particular topological feature in a particular homology group appears
and the $r$ at which the same topological feature disappears as a (birth, death)
pair, which exists in $\R^2$. We can plot these points on a graph to create a
\emph{persistence diagram}, as exemplified in Figure \ref{fig:example_pd}.

While the persistence diagram allows us to visualize the topological features of
a particular sample, it does not extend to separable Banach spaces, and
therefore does not lend itself to analysis through the lens of random
variables\cite{bubenikStatisticalTopologicalData2015}. However, we can transform our
persistence diagram into a \emph{persistence landscape}, which is a function and
does exist in a separable Banach space. We construct persistence landscapes by
drawing an isosceles triangle centered on the points of our persistence
diagrams. When intersections occur the highest function is defined as the
persistence landscape.

\subsection{Permutation Testing}\label{sec:permutation_test}

While we now have usable TDA features in the form of persistence landscapes, we
need to identify a test which can indicate whether our TDA features from OSA
positive pediatric patients is different than our TDA features from OSA negative
patients. To this end, we leverage a permutation test, as described
in~\cite{robinsonHypothesisTestingTopological2017}. In this methodology, our
null hypothesis is that there is no impact of OSA on the brain
connectivity network. To perform the test, we implement the following algorithm:

\begin{enumerate}
  \item Calculate the average OSA positive persistence landscapes
    $\{\lambda_1^{(1)}, \cdots, \lambda_{n_1}^{(1)}\}$ and the average
    OSA negative persistence landscapes $\{\lambda_1^{(2)}, \cdots,
    \lambda_{n_2}^{(2)}\}$.
  \item Calculate the difference between the average OSA positive persistence
    landscape and the average OSA negative persistence landscape, and find the
    most significant value of this difference landscape $\hat F_T$.
  \item Permute the labels of our two groups such that we have $n_1$ randomly
    chosen positively labeled persistence landscapes and the remaining $n_2$
    landscapes are negatively labeled.
  \item Calculate the difference between the permuted positive group and the
    permuted negative group, and find the most significant value of this
    difference landscapes $\hat F_p$.
  \item If $\hat F_p \geq \hat F_T$, count a significant permutation.
  \item Repeat steps 3 and 4 $B$ times. For this analysis we repeat these steps
    $B=1000$ times.
  \item Compute the $p$-value as the number of significant permutations $S$
    divided by the total number of permutations $B$.
\end{enumerate}

\noindent We specifically perform this permutation test on the 0th homology
group ($H_0$) persistence landscapes generated as described in
Section~\ref{sec:tda_methods}.

\section{RESULTS}

\begin{figure}[thpb]
  \framebox{\parbox{3in}{
  \centering
  \includegraphics[scale=0.4]{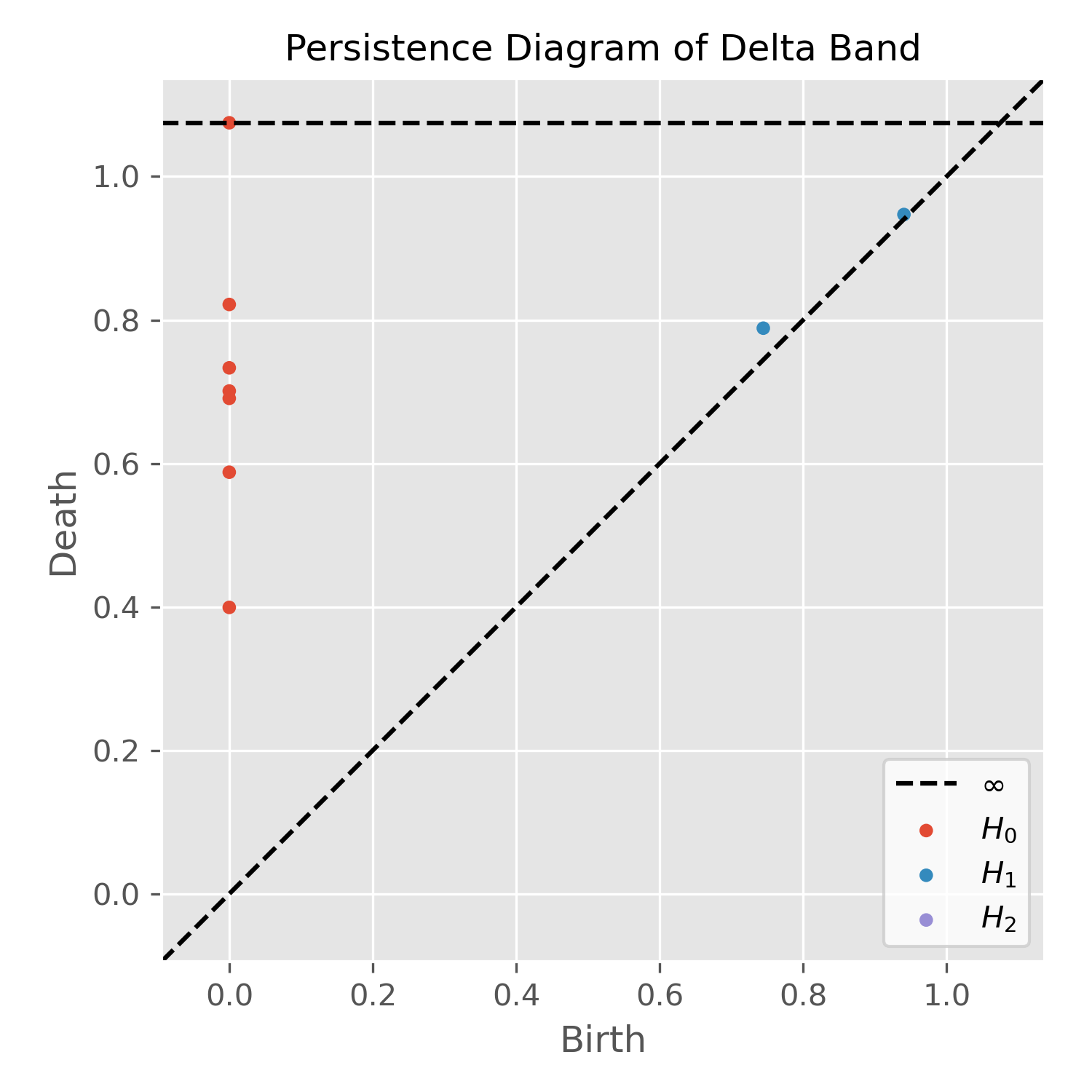}
  \caption{Example Persistence Diagram for Delta Band of an OSA Positive
  Patient}
  \label{fig:example_pd}
  }}
\end{figure}

\begin{figure}[thpb]
  \framebox{\parbox{3in}{
  \centering
  \includegraphics[scale=0.4]{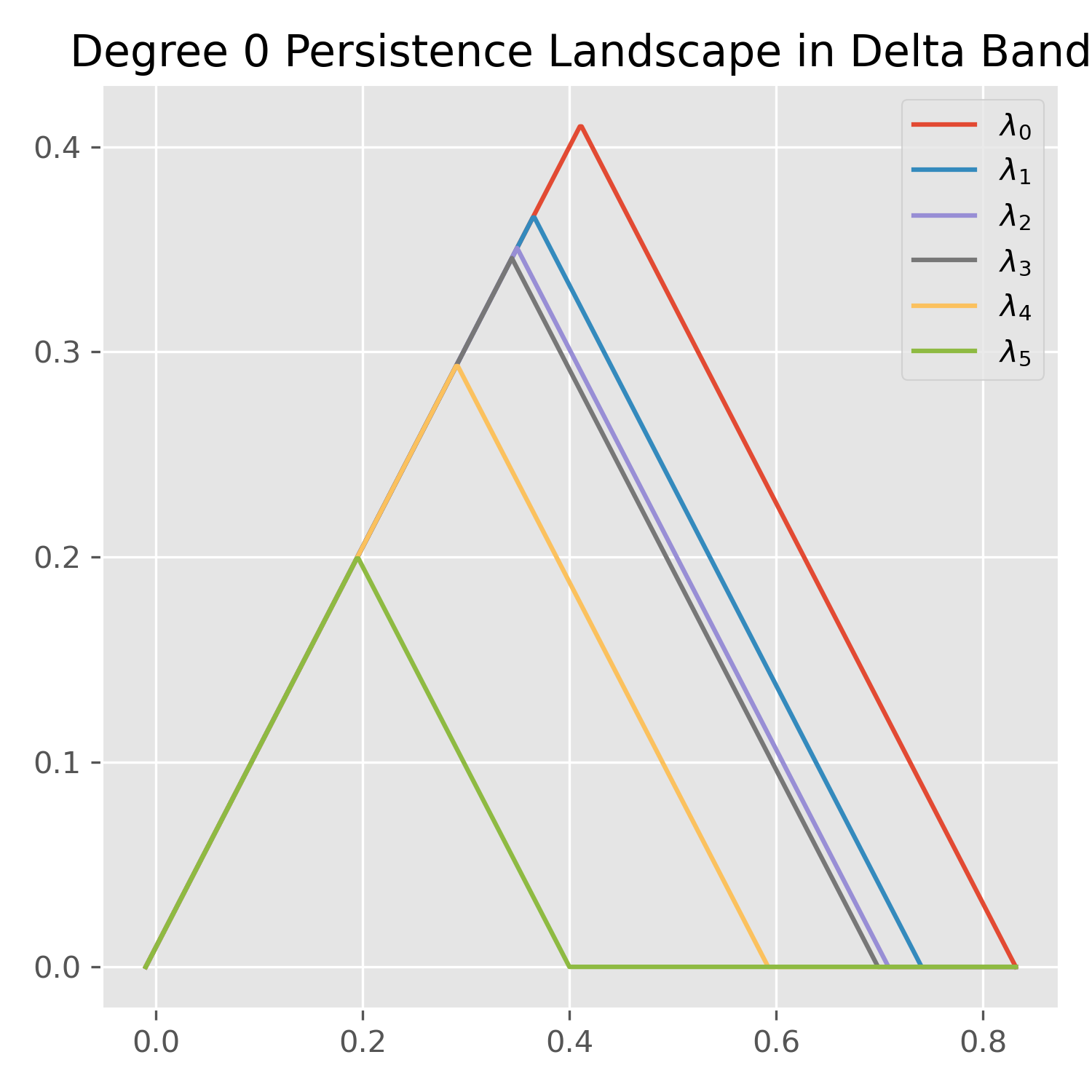}
  \caption{Example Degree 0 Persistence Landscape for Delta Band of an OSA
  Positive Patient}
  \label{fig:example_pl}
  }}
\end{figure}

We show example persistence diagrams and persistence landscapes for an OSA
positive patient in Figure~\ref{fig:example_pd} and Figure~\ref{fig:example_pl}.

We present $p$-values for our TDA permutation tests discussed in
Section~\ref{sec:permutation_test} in Table~\ref{tab:sleep_state_p_value_pl}.
We consider $p$-values with $p < 0.05$ to be significant. We separate our
distance matrices into signals which occur during NREM1 sleep, NREM2 sleep,
NREM3 sleep, and REM sleep.

\begin{table}[H]
	\centering
	\caption{$p$-values for Sleep States using TDA}
	\begin{tabular}{| c | c | c | c | c |}
		\hline
		\textbf{EEG Band} & NREM1  & NREM2  & NREM3  & REM  \\
		\hline
		Delta Band        & 0.000 & 0.001 & 0.000 & 0.040 \\
		\hline
		Theta Band        & 0.000 & 0.000 & 0.000 & 0.000 \\
		\hline
		Alpha Band        & 0.000 & 0.000 & 0.000 & 0.000 \\
		\hline
		Beta Band         & 0.000 & 0.000 & 0.000 & 0.000 \\
		\hline
		Gamma Band        & 0.000 & 0.000 & 0.000 & 0.000 \\
		\hline
	\end{tabular}
  \label{tab:sleep_state_p_value_pl}
\end{table}

\section{DISCUSSION}

Our results indicate that TDA can identify statistically significant differences
between OSA positive and OSA negative cohorts. Statistical significance is
achieved in all bands for all phases of sleep. These results are a strong
indication that EEG signals alone can be used to help diagnose OSA in pediatric
patients, allowing for less invasive diagnosis of OSA.

\section{CONCLUSIONS}

In this work, we show that the brain dependence networks of OSA positive
pediatric patients show a statistically significant difference than the brain
dependence networks of OSA negative pediatric patients when analyzed using
topological data analysis. Our work lays the foundation for using TDA for
cohort-level OSA classification, and takes the first step towards being able to
make an OSA diagnosis with a single EEG. Future work is needed to achieve
clinical utility and to determine how the brain dynamics TDA identifies relate
to respiration measures more typically used for OSA diagnosis.

\addtolength{\textheight}{-12cm}   



%

\section*{ACKNOWLEDGMENT}

NCH Sleep DataBank was supported by the National Institute of Biomedical Imaging
and Bioengineering of the National Institutes of Health under Award Number
R01EB025018. The National Sleep Research Resource was supported by the U.S.
National Institutes of Health, National Heart Lung and Blood Institute (R24
HL114473, 75N92019R002).


\bibliographystyle{plain}
\bibliography{smanjunath_arxiv}

\end{document}